\documentclass[12pt]{article}
\usepackage{exscale,relsize}
\usepackage{amsmath}
\usepackage[sort,nocompress]{cite}
\usepackage{amsfonts}
\usepackage[hidelinks]{hyperref}
\usepackage{amssymb}
\usepackage{calc}
\usepackage{theorem}
\usepackage{pifont}      
\usepackage{array}
\usepackage{color}

\usepackage{graphicx} 
\usepackage{float} 
\usepackage{subcaption} 

\newcommand{\Ball}{\mathbb{B}}

\newcommand{\N}{\mathbb{N}}   

\newcommand{\dist}{\ensuremath{\operatorname{dist}}} 

\newcommand{\norm}[1]{\left\lVert#1\right\rVert} 
\newcommand{\ip}[2]{\langle#1,#2\rangle} 

\newcommand{\R}{\ensuremath{\mathbb R}}
\newcommand{\Rn}{\ensuremath{\mathbb R^n}}

\newcommand{\dom}{\ensuremath{\operatorname{dom}}}

\newcommand{\inte}{\ensuremath{\operatorname{int}}}

\providecommand{\BB}[2]{\mathbb{B}(#1;#2)}

\newtheorem{theorem}{Theorem}[section]
\newtheorem{lemma}[theorem]{Lemma}
\newtheorem{fact}[theorem]{Fact}
\newtheorem{corollary}[theorem]{Corollary}
\newtheorem{proposition}[theorem]{Proposition}
\newtheorem{defn}[theorem]{Definition}

\theoremstyle{plain}{\theorembodyfont{\rmfamily}
}
\theoremstyle{plain}{\theorembodyfont{\rmfamily}
}
\theoremstyle{plain}{\theorembodyfont{\rmfamily}
}
\theoremstyle{plain}{\theorembodyfont{\rmfamily}
\newtheorem{example}[theorem]{Example}}
\theoremstyle{plain}{\theorembodyfont{\rmfamily}
\newtheorem{remark}[theorem]{Remark}}
\theoremstyle{plain}{\theorembodyfont{\rmfamily}
}
\def\proof{\noindent{\it Proof}. \ignorespaces}
\def\endproof{\ensuremath{\quad \hfill \blacksquare}}

\newcommand{\pluss}{{\hskip1pt \raise1pt\vbox{\hrule width6pt \vskip1pt
\hrule width6pt}\kern-4pt{\lower1pt\hbox{\vrule height6pt \kern1pt\vrule
height6pt}}\hskip5pt}}

\begin{document}
\title{{\fontfamily{ptm}\selectfont Calculus rules of the generalized concave Kurdyka-\L ojasiewicz property}}
\author{Xianfu Wang\footnote{Mathematics, University of British Columbia, Kelowna, B.C. V1V 1V7, Canada. E-mail: shawn.wang@ubc.ca.} and Ziyuan Wang\footnote{Mathematics, University of British Columbia, Kelowna, B.C. V1V 1V7, Canada. E-mail: ziyuan.wang@alumni.ubc.ca.}}
\date{August 26, 2021}
\maketitle
\begin{abstract}\noindent In this paper, we propose several calculus rules for the generalized concave Kurdyka-\L ojasiewicz~(KL) property, which generalize Li and Pong's results for KL exponents. The optimal concave desingularizing function has various forms and may be nondifferentiable. Our calculus rules do not assume desingularizing functions to have any specific form nor differentiable, while the known results do. Several examples are also given to show that our calculus rules are applicable to a broader class of functions than the known ones.
\end{abstract}


\section{Introduction}\label{Sec: intro}
Throughout this paper,
\[\Rn\text{ is the standard Euclidean space}\]
equipped with inner product $\ip{x}{y}=x^Ty$ and the Euclidean norm $\norm{x}=\sqrt{\ip{x}{x}}$ for every~$x,y\in\Rn$. Denote by $\N$ the set of positive natural numbers, i.e., $\N=\{1,2,3,\ldots\}$. The open ball centered at $\bar{x}$ with radius $r$ is denoted by~$\BB{\bar{x}}{r}$. The distance function of a subset $K\subseteq \Rn$ is $\dist(\cdot,K):\Rn\rightarrow\overline{\R}=(-\infty,\infty]$,
\[x\mapsto\dist(x,K)=\inf\{\norm{x-y}:y\in K\},\]
where $\dist(x,K)\equiv\infty$ if $K=\varnothing $. For $f:\Rn\to\overline{\R}$ and $r_1,r_2\in[-\infty,\infty]$, we set $[r_1<f<r_2]=\{x\in\Rn:r_1<f(x)<r_2\}$. For $\eta\in(0,\infty]$, we denote by $\Phi_\eta$ the class of functions $\varphi:[0,\eta)\rightarrow\R_+$ satisfying the following conditions: (i) $\varphi(t)$ is right-continuous at $t=0$ with $\varphi(0)=0$; (ii) $\varphi$ is strictly increasing on $[0,\eta)$. Recall that the left derivative of $\varphi:[0,\infty)\rightarrow\R$ at $t$ is
\[\varphi_-^\prime(t)=\lim_{s\rightarrow t^-}\frac{\varphi(s)-\varphi(t)}{s-t}.\]
Before stating the goal of this paper, let us recall the pointwise version of the concave Kurdyka-\L ojasiewicz~(KL) property. Let $\mathcal{K}_\eta$ denote all functions in $\Phi_\eta$ that are continuously differentiable on $(0,\eta)$.
\begin{defn}\label{Def:KL property}
	Let $f:\mathbb{R}^n\rightarrow\overline{\mathbb{R}}$ be proper and lsc. We say $f$ has the concave\footnote{Unlike many published articles, we use the adjective ``concave" because the concavity of desingularizing function is an additional property useful for algorithmic applications and was not assumed in the seminal work on KL property, see, e.g.,~\cite{bolte2010survey}. We appreciate the anonymous referee of the companion paper~\cite{wang2020} for suggesting this terminology.} KL property at $\bar{x}\in\dom\partial f$, if there exist neighborhood $U\ni\bar{x}$, $\eta\in(0,\infty]$ and a concave function $\varphi\in\mathcal{K}_\eta$ such that for all $x\in U\cap[0<f-f(\bar{x})<\eta]$,
	\begin{equation}\label{Inequality: KL inequality}
		\varphi^\prime\big(f(x)-f(\bar{x})\big)\cdot\dist\big(0,\partial f(x)\big)\geq1,
	\end{equation}		
	where $\partial f(x)$ denotes the limiting subdifferential of $f$ at $x$, see Definition~\ref{Defn:limiting subdifferential}. The function $\varphi$ is called a concave \textit{desingularizing function} of $f$ at $\bar{x}$ with respect to $U$ and $\eta$. We say $f$ is a concave KL function if it has the concave KL property at every $\bar{x}\in\dom\partial f$.
\end{defn}
The concave KL property is instrumental in the convergence analysis
of many proximal-type algorithms, see,~e.g.,~\cite{Attouch2009,attouch2010proximal,attouch2013convergence,
Bolte2014,Ipiano2016,TKP2017extra,Banert2019,yu2019deducing,yu2020convergence,boct2020extrapolated} and the references therein; see also~\cite{bolte2010survey,Lojas1963,Kur98} for seminal theoretical work on this area. Convergence rates of such algorithms are usually determined by the KL exponent $\theta\in[0,1)$ when desingularizing functions have the \L ojasiewicz form~$\varphi(t)=c\cdot t^{1-\theta}$ for some $c>0$. Calculating the KL exponent is usually challenging, however there are a few positive results.
Li and Pong~\cite{li2018calculus} recently developed several important
calculus rules for the KL exponent, which facilitates estimating the KL exponent for some structured optimization problems. Under suitable assumptions, Wu, Pan and Bi~\cite{wu2021kurdyka} studied the KL exponent for problems involving sums of the zero norm and nonconvex functions, see also~\cite{yu2019deducing,liu2019quadratic} for other pleasing progress on this line of research.

\emph{Our goal in this paper is to develop calculus rules
that facilitate finding desingularizing functions that are neither differentiable nor have any specific forms}. In a companion paper~\cite{wang2020}, we introduced the generalized concave KL property, which generalizes the concave KL property, see Definition~\ref{Def: g-KL} and Remark~\ref{rem: g-Kl}. It turns out that this generalized framework allows us to capture the smallest concave desingularizing function through its associated exact modulus. In particular, the optimal concave desingularizing function may neither have the \L ojasiewicz form nor be differentiable, see Example~\ref{ex: nondifferentiable modulus} and Proposition~\ref{Prop: the usual form doesn't always work}, which naturally motivates this work. Moreover, our goal is interesting in its own right since most published results emphasize on calculus rules for desingularizing functions with the  \L ojasiewicz form, while calculus rules for general desingularizing functions received little attention and require a more sophisticated construction. Our~\textbf{main contributions} are listed below:
\begin{itemize}

  \item Theorem~\ref{thm: multiple sum rule} states that the sum $f=\sum_{i=1}^mf_i$ admits a concave desingularizing function $\varphi(t)=\frac{1}{\alpha}\int_0^t\max_{1\leq i\leq m}(\varphi_i)^\prime_-\left(\frac{s}{m}\right)ds$ at $\bar{x}$ for some $\alpha>0$, provided that $\varphi_i$ is a concave desingularizing function for $f_i$ at $\bar{x}$ for each $i$ and a regularity condition~(\ref{linear regularity of m functions}) is satisfied.

  \item Theorem~\ref{Thm:minmum of g-KL} shows that $\varphi(t)=\int_0^t\max_{i\in I(\bar{x})}(\varphi_i)_-^\prime(s)ds$ is a concave desingularizing function of
	 $f(x)=\min_{1\leq i\leq m}f_i(x)$ at $\bar{x}$, given that $f_i$ has a concave desingularizing function $\varphi_i$ for each $i$, where $I(\bar{x})=\{i:f_i(\bar{x})=f(\bar{x}), 1\leq i\leq m\}$.
  \item Theorem~\ref{Thm:separable sum} concerns the generalized concave KL property of separable sums. The function $\varphi(t)=\int_0^t\max_{1\leq i\leq m}\left[(\varphi_i)^\prime_-\left(\frac{s}{m}\right)\right]ds$ is a concave desingularizing function of $f(x)=\sum_{i=1}^mf_i(x_i)$ at $\bar{x}=(\bar{x}_1,\ldots,\bar{x}_m)$, provided that $f_i$ has a concave desingularizing function $\varphi_i(t)$ at $\bar{x}_i$ for each $i$.
  \item A composition rule is given in Theorem~\ref{Thm: Composition rule}. For $f(x)=(g\circ F)(x)$, where $F$ is a smooth map and $g$ has a concave desingularizing function $\varphi$ at $\bar{y}=F(\bar{x})$, there exists $r>0$ such that $f$ has a concave desingularizing function $\varphi(t)/r$ at $\bar{x}$.
\end{itemize}
 Theorems~\ref{thm: multiple sum rule},~\ref{Thm:minmum of g-KL},~\ref{Thm:separable sum} and~\ref{Thm: Composition rule}
 significantly generalize corresponding calculus rules by Li and Pong, see~\cite{li2018calculus} and its early preprint~\cite{li2016calculus}. Examples are given to demonstrate that they are applicable to a broader class of functions, see Example~\ref{ex:min rule fails}.


The structure of this paper is as the following: In Section~\ref{Section:Preliminaries}, we collect background knowledge about variational analysis. Our main results and examples are presented in Section~\ref{Section:calculus rules}. We end this paper in Section~\ref{Section:Conclusion} with a concluding remark and discussion about future work. Pleasant properties of the generalized KL property and its associated exact modulus, as well as supplementary lemmas, can be found in the Appendix.

\section{Preliminaries and facts}\label{Section:Preliminaries}
In this paper, we will use the following generalized subgradients, see, e.g., \cite{mor2018variational,rockwets}.
\begin{defn}\label{Defn:limiting subdifferential}
	Let $f:\Rn\rightarrow\overline{\R}$ be a proper function. We say that
	\begin{itemize}
		\item[(i)] $v\in\Rn$ is a \textit{Fr\'echet subgradient} of $f$ at $\bar{x}\in\dom f$, denoted by $v\in\hat{\partial}f(\bar{x})$, if for every $x\in\dom f$,
		\begin{equation}\label{Formula:frechet subgradient inequality}
		f(x)\geq f(\bar{x})+\ip{v}{x-\bar{x}}+o(\norm{x-\bar{x}}).
		\end{equation}
		\item[(ii)] $v\in\Rn$ is a \textit{limiting subgradient} of $f$ at $\bar{x}\in\dom f$, denoted by $v\in\partial f(\bar{x})$, if
		\begin{equation}\label{Formula:limiting subgraident definition}
		v\in\{u\in\Rn:\exists x_k\xrightarrow[]{f}\bar{x},\exists u_k\in\hat{\partial}f(x_k),u_k\rightarrow u\},
		\end{equation}
where $x_k\xrightarrow[]{f}\bar{x}$ denotes $x_k\to\bar{x}$ and $f(x_k)\to f(\bar{x})$ as $k\to\infty$. Moreover, we set $\dom\partial f=\{x\in\Rn:\partial f(x)\neq\varnothing \}$. We say $\bar{x}$ is a stationary point, if $0\in\partial f(\bar{x})$.
	\end{itemize}
\end{defn}


The following subdifferential calculus rules will used in the sequel.

\begin{fact}\label{Fact: subdifferential of min} \emph{\cite[Proposition 4.9]{mor2018variational}} Suppose that $f(x)=\min_{1\leq i\leq m} f_i(x)$, where $f_i:\Rn\rightarrow\overline{\R}$ is lsc and proper function. Let $x\in\bigcap_{i=1}^m
\dom\partial f_i$ and $I(x)=\{i:f(x)=f_i(x),1\leq i\leq m\}$. Then
	\[\partial f(x)\subseteq\bigcup_{i\in I(x)}\partial f_i(x).\]
	\end{fact}

\begin{fact}\label{Fact:subd for separable sum}\emph{\cite[Proposition 10.5]{rockwets}}
Suppose that $f(x)=\sum_{i=1}^m f_i(x_i)$, where $f_i:\R^{n_i}\rightarrow\overline{\R}$ is proper and lsc for each $i\in\{1,\ldots,m\}$. Then we have \[\partial f(x)=\prod_{i=1}^{m}\partial f_i(x_i).\]
\end{fact}

\begin{fact}\label{Fact: chain rule}\emph{\cite[Excercise 10.7]{rockwets}} Suppose that $f(x)=g(F(x))$ for a proper lsc function $g:\R^m\rightarrow\overline{\R}$ and a smooth map $F:\Rn\rightarrow\R^m$. Let $\bar{x}$ be a point where $f$ is
finite and the Jacobian $\nabla F(\bar{x})$ is surjective. Then
	\begin{equation}
\hat{\partial}f(\bar{x})=\nabla F(\bar{x})^*\hat{\partial}g(F(\bar{x}))\text{ and }\partial f(\bar{x})=\nabla F(\bar{x})^*\partial g(F(\bar{x})).
	\end{equation}
	\end{fact}
\begin{lemma}\label{Lemma: properties of convex functions on the line}\emph{\cite[Theorems 4.42--4.43]{Stromberg}} Let $I\subseteq\R$ be an open interval and let $\varphi:I\rightarrow\R$ be convex. Then
the following hold:

(i) The side derivatives $\varphi_-^\prime(t)$ and $\varphi_+^\prime(t)$ are finite at every $t\in I$. Moreover, $\varphi_-^\prime(t)$ and $\varphi_+^\prime(t)$ are increasing.

(ii) $\varphi$ is differentiable almost everywhere on $I$.

(iii) Let $t\in I$. Then for every $s\in I$, $\varphi(s)-\varphi(t)\geq \varphi_-^\prime(t)\cdot(s-t)$.
\end{lemma}



\section{Calculus rules of the generalized Kurdyka-\L ojasiewicz property}\label{Section:calculus rules}

We shall develop desired calculus rules within the framework of the generalized concave KL property introduced in~\cite{wang2020}, which is a generalization of the concave KL property (recall Definition~\ref{Def:KL property}). See Appendix~\ref{AppendixA} for their pleasant properties.

\begin{defn}\label{Def: g-KL} Let $f:\Rn\rightarrow\overline{\mathbb{R}}$ be proper lsc and let $\bar{x}\in\dom\partial f$.
	
(i)	We say that $f$ has the \textit{generalized concave Kurdyka-\L ojasiewicz property} at $\bar{x}\in\dom\partial f$, if there exist neighborhood $U\ni\bar{x}$, $\eta\in(0,\infty]$ and a concave $\varphi\in\Phi_\eta$, such that for all $x\in U\cap[0<f-f(\bar{x})<\eta]$,
	\begin{equation}\label{g-KL inequality}
		\varphi^\prime_-\big(f(x)-f(\bar{x})\big)\cdot\dist\big(0,\partial f(x)\big)\geq1.
	\end{equation}

(ii) Define $h:(0,\eta)\rightarrow\R$ by
\[h(s)=\sup\big\{\dist^{-1}\big(0,\partial f(x)\big): x\in U\cap[0<f-f(\bar{x})<\eta],s\leq f(x)-f(\bar{x})\big\}.\]
Suppose that $h(s)<\infty$ for $s\in(0,\eta)$. The exact modulus of the generalized concave KL property of $f$ at $\bar{x}$ with respect to $U$ and $\eta$ is defined by the function
$\tilde{\varphi}:[0,\eta)\rightarrow\R_+$,	\begin{align}\label{g-KL modulus}t\mapsto\int_0^th(s)ds,~\forall t\in(0,\eta),	\end{align}
and $\tilde{\varphi}(0)=0$. If $U\cap[0<f-f(\bar{x})<\eta]=\varnothing $ for given $U\ni\bar{x}$ and $\eta>0$, then we set the exact modulus with respect to $U$ and $\eta$ to be $\tilde{\varphi}(t)\equiv0$.
\end{defn}
\begin{remark}\label{rem: g-Kl} The left derivative $t\mapsto\varphi_-^\prime(t)$ is well defined on $(0,\eta)$, see Lemma~\ref{Lemma: properties of convex functions on the line}(i). When $\varphi$ is continuously differentiable on $(0,\eta)$, it is easy to see that the generalized concave KL property reduces to the concave KL property.
\end{remark}

We begin our main results with a sum rule, which generalizes~\cite[Theorem 3.4]{li2016calculus} under a weaker condition.
\begin{theorem}\label{thm: multiple sum rule}Let $f_i:\Rn\to\overline{\R}$ be a proper continuous function for each $i\in\{1,\ldots,m\}$ with $\cap_{i=1}^m\inte\dom\partial f_i\neq\varnothing$, and let $f=\sum_{i=1}^mf_i$. Suppose that at most  one of $f_i$ is not locally Lipschitz. Pick $\bar{x}\in\cap_{i=1}^m\inte\dom\partial f_i$. Assume that for each $i$, $f_i$ has the generalized concave KL property at $\bar{x}$ with respect to $U_i=\BB{\bar{x}}{\varepsilon_i}$ for some $\varepsilon_i>0$, $\eta_i\in(0,\infty]$ and $\varphi_i\in\Phi_{\eta_i}$. Suppose that there exist $\varepsilon_0>0$ and $\alpha>0$ such that for every $x_i\in\BB{\bar{x}}{\varepsilon_0}$
	\begin{equation}\label{linear regularity of m functions}
		\norm{\sum_{i=1}^mu_i}\geq\alpha\sum_{i=1}^m\norm{u_i},~\forall u_i\in\partial f_i(x_i).
	\end{equation}
	Set $\varepsilon=\min_{0\leq i\leq m}\varepsilon_i$ and $\eta=\min_{1\leq i\leq m}\eta_i$. Define $\varphi:[0,\eta)\to\R$ by
	\[\varphi(t)=\frac{1}{\alpha}\int_0^t\max_{1\leq i\leq m}(\varphi_i)^\prime_-\left(\frac{s}{m}\right)ds,\forall t\in(0,\eta),\]
	and $\varphi(0)=0$. If each $\varphi_i$ is strictly concave, then the sum $f$ has the generalized concave KL property at $\bar{x}$ with respect to $U=\BB{\bar{x}}{\varepsilon}$, $\eta$ and $\varphi$.
\end{theorem}
\proof By Lemma~\ref{Lemma: integral of maxmum of left-d is well-defined}, one concludes that $\varphi$ is well-defined and right-continuous at $0$. For convenience, we write $\phi_i=(\varphi_i)_-^\prime$ for every $i$. Invoking Lemma~\ref{Lemma:integration of increasing function is convex}, one concludes that $\varphi$ belongs to $\Phi_\eta$ and
\begin{equation}\label{a}
	\varphi_-^\prime(t)\geq\frac{1}{\alpha}\max_{1\leq i\leq m}\phi_i\left(\frac{t}{m}\right)\geq\frac{1}{\alpha}\phi_i\left(\frac{t}{m}\right),\forall 1\leq i\leq m.
\end{equation}
Note $\varphi_i$ is strictly concave. Then Lemma~\ref{left derivative inequality} ensures that the function $\phi_i$ is strictly decreasing and hence invertible. By assumption, assume without loss of generality that $f_m$ is not locally Lipschitz. Picking $x\in U$, one sees easily that the set $\sum_{i=1}^{m-1}\partial f_i(x) $ is compact while $\partial f_m(x)$ is closed. Hence $\sum_{i=1}^m\partial f_i(x)$ is closed and there exists $u_i\in\partial f_i(x)$ for each $i$ such that $\dist\left(0,\sum_{i=1}^m\partial f_i(x)\right)=\norm{\sum_{i=1}^mu_i}$. Then~(\ref{linear regularity of m functions}) implies that
\begin{align}\label{x}
	\dist(0,\partial f(x))\geq \dist\left(0,\sum_{i=1}^m\partial f_i(x)\right)=\norm{\sum_{i=1}^mu_i}\geq \alpha\sum_{i=1}^m\norm{u_i}\geq \alpha\dist(0,\partial f_i(x)),
\end{align}
where the first inequality holds because of our assumption and the subdifferential sum rule, see, e.g.,~\cite[Corollary 10.9]{rockwets}.

Pick $x\in U$ and assume without loss of generality that $f_i(x)-f_i(\bar{x})<\eta$ for every $i\in\{1,\ldots,m\}$. We claim that for each $i\in\{1,\ldots,m\}$,
\begin{equation}\label{xxxxx}
	f_i(x)-f_i(\bar{x})\leq\phi_i^{-1}\left(\frac{\alpha}{\dist(0,\partial f(x))}\right),
\end{equation}
which holds trivially when $f_i(x)\leq f_i(\bar{x})$. For $x\in U$ with $f_i(x)>f_i(\bar{x})$, one has $x\in U_i\cap[0<f-f_i(\bar{x})<\eta_i]$. By the assumption that $f_i$ has the generalized concave KL property with respect to $U_i$, $\eta_i$ and $\varphi_i$, we have
\begin{equation*}
	\phi_i\left(f_i(x)-f_i(\bar{x})\right)\geq \frac{1}{\dist(0,\partial f_i(x))}\geq\frac{\alpha}{\dist(0,\partial f(x))},
\end{equation*}
where the last inequality follows from~(\ref{x}), from which the claimed inequality~(\ref{xxxxx}) readily follows. Now, for $x\in U$, let $i^*=i(x)$ be the index such that for every $i$ \[\phi_{i^*}^{-1}\left(\frac{\alpha}{\dist(0,\partial
	f(x))}\right)\geq\phi_i^{-1}\left(\frac{\alpha}{\dist(0,\partial f(x))}\right).\]
Then we have
\begin{equation}\label{z}
	f(x)-f(\bar{x})=\sum_{i=1}^mf_i(x)-f_i(\bar{x})\leq \sum_{i=1}^m\phi_i^{-1}\left(\frac{\alpha}{\dist(0,\partial f(x))}\right)\leq m\phi_{i^*}^{-1}\left(\frac{\alpha}{\dist(0,\partial f(x))} \right),
\end{equation}
where the first inequality is implied by~(\ref{xxxxx}).

Finally, we prove the desired statement. For $x\in\BB{\bar{x}}{\varepsilon}\cap[0<f-f(\bar{x})<\eta]$, we have
\begin{align*}
	\varphi_-^\prime\left(f(x)-f(\bar{x})\right)\dist(0,\partial f(x))&\geq\frac{1}{\alpha}\phi_{i^*}\left(\frac{f(x)-f(\bar{x})}{m}\right)\dist(0,\partial f(x))\\
	&\geq\frac{1}{\alpha}\phi_{i^*}\left(\phi_{i^*}^{-1}\left(\frac{\alpha}{\dist(0,\partial f(x))} \right)\right)\dist(0,\partial f(x))\\
	&=1,
\end{align*}
where the first inequality is implied by~(\ref{a}) and the second holds because of~(\ref{z}).\endproof

\begin{remark}\label{rem: sum rule} (i) Condition~(\ref{linear regularity of m functions}) is weaker than the Li-Pong disjoint condition~\cite[Condition (12)]{li2016calculus} for a sum rule of two functions. Considering $f=f_1+f_2$ with associated desingularizing functions having the form $\varphi_i(t)=c_i\cdot t^{1-\theta_i}$ for some $c_i>0$ and $\theta_i\in[0,1)$, Li and Pong provided a sum rule~\cite[Theorem 3.4]{li2016calculus} assuming the following:
	\begin{equation}\label{li-pong disjoint condition}	W_{f_1}(\bar{x})\cap(-W_{f_2}(\bar{x}))=\varnothing,\end{equation}where $$W_{f_i}(\bar{x})=\bigg\{\lim_{k\to\infty}w_k:w_k=\frac{u_k}{\norm{u_k}}, u_k\in\partial f_i(x_k)\backslash\{0\}, x_k\to\bar{x}\bigg\},$$
which implies~(\ref{linear regularity of m functions}) with $m=2$ by the first line of~\cite[Inequality (16)]{li2016calculus}. The converse, however, fails even on the real line. For instance, consider $f_1(x)=x$ and $f_2(x)=-\frac{1}{2}x$. On one hand, we have $|f_1^\prime(x)+f_2^\prime(x)|=\frac{1}{3}\left(|f_1^\prime(x)|+|f_2^\prime(x)|\right)$ for every $x\in\R$. On the other hand, it is easy to see that $W_{f_1}(0)=-W_{f_2}(0)=\{1\}$, which means that~(\ref{li-pong disjoint condition}) fails. Nevertheless, these two conditions are equivalent when $f_i=\delta_{C_i}$, where $C_i$ are nonempty closed sets, see, e.g.,~\cite[Proposition 2.4]{dao2019linear}.
	
(ii) For $i=\{1,\ldots,m\}$, if $\eta_i<m$ and  $\varphi_i(t)=t^{1-\theta_i}/(1-\theta_i)$ for some $\theta_i\in[0,1)$, then the desingularizing function given by Theorem~\ref{thm: multiple sum rule} reduces to \[\varphi(t)=\frac{1}{\alpha}\int_0^t\max_{1\leq i\leq m}\left(\frac{s}{m} \right)^{-\theta_i}ds=\frac{1}{\alpha}\int_0^t\left(\frac{s}{m}\right)^{-\theta} ds=\frac{m^{\theta}}{(1-\theta)\alpha}t^{1-\theta},\forall t\in(0,\eta)\subseteq(0,m),\]
where $\theta=\max_{1\leq i\leq m}\theta_i$ and the second equality holds because $s/m<1$.
\end{remark}

Taking Theorem~\ref{thm: multiple sum rule} and Remark~\ref{rem: sum rule}(ii) into account, one can immediately obtain the following corollary, which states the same conclusion as in~\cite[Theorem 3.4]{li2016calculus} but under a weaker condition.

\begin{corollary} Let $f_1,f_2:\Rn\to\overline{\R}$ be proper continuous functions. Assume that $f_1$ is locally Lipschitz. Pick $\bar{x}\in\inte\dom\partial f_1\cap\inte\dom\partial f_2$ and suppose that for each $i\in\{1,2\}$ the function $f_i$ has the KL exponent $\theta_i\in(0,1)$ at $\bar{x}$. Assume further that~(\ref{linear regularity of m functions}) holds for some $\alpha>0$. Then the sum $f=f_1+f_2$ has the KL exponent $\theta=\max(\theta_1,\theta_2)$.
\end{corollary}

The following theorem is a generalization of~\cite[Theorem 3.1]{li2018calculus}, whose proof we follow. The novelty of our result is that we put no restriction on the differentiability and form of desingularizing function.

\begin{theorem}\label{Thm:minmum of g-KL} Suppose that $f(x)=\min_{1\leq i\leq m}f_i(x)$ is continuous on $\dom\partial f$, where $f_i:\Rn\rightarrow\overline{\R}$ is proper and lsc function for each $i\in\{1,\ldots,m\}$. Let $\bar{x}\in\dom\partial f$ and suppose that $\bar{x}\in\bigcap_{i\in I(\bar{x})}\dom\partial f_i$, where $I(\bar{x})=\{i:f_i(\bar{x})=f(\bar{x})\}$. Assume further that for every $i\in I(\bar{x})$, the function $f_i$ has the generalized concave KL property at $\bar{x}$ with respect to $U_i=\BB{\bar{x}}{\varepsilon_i}$ for some $\varepsilon_i>0$, $\eta_i>0$ and $\varphi_i\in\Phi_{\eta_i}$. Set $\eta=\min_{1\leq i\leq m}\eta_i$ and define $\varphi:[0,\eta)\rightarrow\overline{\R}$ by
	\[\varphi(t)=\int_0^t\max_{i\in I(\bar{x})}(\varphi_i)_-^\prime(s)ds,\forall 0<t<\eta,\]
and $\varphi(0)=0$.	Then there exists $\varepsilon>0$ such that $f$ has the generalized concave KL property at $\bar{x}$ with respect to $U=\BB{\bar{x}}{\varepsilon}$, $\eta$ and $\varphi$.
	\end{theorem}
\proof Note that $\varphi_i(t)$ is finite and right-continuous at $0$ for each $i$. Hence Lemma~\ref{Lemma: integral of maxmum of left-d is well-defined} ensures that $\varphi(t)$ is well-defined and continuous at $0$. On the other hand, for every $i$, we have $\varphi_i(t)=\int_0^t(\varphi_i)_-^\prime(s)ds$, where $(\varphi_i)_-^\prime(t)$ is a decreasing function. Then $t\mapsto\max(\varphi_i)_-^\prime(t)$ is a decreasing function, which by Lemma~\ref{Lemma:integration of increasing function is convex} implies that $\varphi\in\Phi_\eta$ and for $t\in(0,\eta)$
	 \begin{equation}\label{dfads} \varphi_-^\prime(t)\geq\max_{i\in I(\bar{x})}(\varphi_i)_-^\prime(t).
	  \end{equation}
Now we work towards the existence of some $\varepsilon_0>0$ such that \begin{equation}\label{Formula: active constraints inclusion}
		I(x)\subseteq I(\bar{x}),\forall x\in\BB{\bar{x}}{\varepsilon_0}.
		\end{equation}
 We claim that there exists $\varepsilon_0>0$ such that \[\min_{i\notin I(\bar{x})}f_i(x)> f(x),\forall x\in\BB{\bar{x}}{\varepsilon_0},\]
 which implies (\ref{Formula: active constraints inclusion}). To see this implication, taking $i_0\in I(x)$ and supposing that $i_0\notin I(\bar{x})$, our claim enforces that $f_{i_0}(x)\geq \min_{i\notin I(\bar{x})} f_i(x)>f(x)=f_{i_0}(x)$ whenever $x\in\BB{\bar{x}}{\varepsilon_0}$, which is absurd. Next we justify the aforementioned claim. Suppose to the contrary that for every $\varepsilon>0$ there exists some $x\in\BB{\bar{x}}{\varepsilon}$ such that
		$\min_{i\notin I(\bar{x})}f_i(x)\leq f(x)$.
		Then there exists a sequence $x_k\rightarrow\bar{x}$ with $$\min_{i\notin I(\bar{x})}f_i(x_k)\leq f(x_k),\forall k\in\N.$$
		Taking $k\rightarrow\infty$, one has by the lower semi-continuity of $x\mapsto\min_{i\notin I(\bar{x})} f_i(x)$ and the continuity of $f(x)$ that
		\[\min_{i\notin I(\bar{x})} f_i(\bar{x})\leq\liminf_{k\rightarrow\infty}\left(\min_{i\notin I(\bar{x})}f_i(x_k)\right)\leq\liminf_{k\rightarrow\infty}f(x_k)=f(\bar{x}),\]
		which is absurd because $\min_{i\notin I(\bar{x})} f_i(\bar{x})> f(\bar{x})$.
		
To show that $\varphi(t)$ is a desingularizing function, inequality~(\ref{Formula:dfasf}) helps, which will be proved in the sequel. By assumption one has for every $i\in I(\bar{x})$ and $x\in\BB{\bar{x}}{\varepsilon_i}\cap[0<f_i-f_i(\bar{x})<\eta_i]$
				\begin{equation}\label{Formula:fklsdjf }
				\dist(0,\partial f_i(x))\geq\frac{1}{(\varphi_i)_-^\prime(f_i(x)-f_i(\bar{x}))}.
				\end{equation}
		Moreover, invoking Fact~\ref{Fact: subdifferential of min}, we have for every $x\in\dom\partial f$
				\[\partial f(x)\subseteq\bigcup_{i\in I(x)}\partial f_i(x).\]
				Take $u\in\partial f(x)$ with $\norm{u}=\dist(0,\partial f(x))$. Then the inclusion above implies that $u\in\partial f_{i_0}(x)$ for some $i_0\in I(x)$, and consequently
				\begin{equation}\label{Formula:dfasf}
				\dist(0,\partial f(x))=\norm{u}\geq\dist(0,\partial f_{i_0}(x))\geq\min_{i\in I(x)}\dist(0,\partial f_i(x)).
				\end{equation}
Set $\varepsilon=\min_{0\leq i\leq m}\varepsilon_i$ and $\eta=\min_{1\leq i\leq m}\eta_i$. Take $x\in\BB{\bar{x}}{\varepsilon}\cap[0<f-f(\bar{x})<\eta]$. Then we have $x\in\BB{\bar{x}}{\varepsilon_i}$ and $0<f(x)-f(\bar{x})=f_i(x)-f_i(\bar{x})<\eta_i$ for every $i\in I(x)\subseteq I(\bar{x})$.
		
Now we are ready to prove that $\varphi$ is a desingularizing function of $f$ at $\bar{x}$. For $x\in\BB{\bar{x}}{\varepsilon}\cap[0<f-f(\bar{x})<\eta]$, we have
		\begin{align*}
		\dist(0,\partial f(x))&\geq\min_{i\in I(x)}\dist(0,\partial f_i(x))\geq\min_{i\in I(x)}\frac{1}{(\varphi_i)_-^\prime(f(x)-f(\bar{x}))}\\
		&\geq\min_{i\in I(\bar{x})}\frac{1}{(\varphi_i)_-^\prime(f(x)-f(\bar{x}))}=\frac{1}{\max_{i\in I(\bar{x})}(\varphi_i)_-^\prime(f(x)-f(\bar{x}))}\\
		&\geq\frac{1}{\varphi_-^\prime(f(x)-f(\bar{x}))},
		\end{align*}
		where the first inequality is implied by~(\ref{Formula:dfasf}); the second one holds because of~(\ref{Formula:fklsdjf }) and the fact that $x\in\BB{\bar{x}}{\varepsilon_i}\cap[0<f_i-f_i(\bar{x})<\eta_i]$ for $i\in I(x)\subseteq I(\bar{x})$; the third one is implied by (\ref{Formula: active constraints inclusion}) and the last inequality holds because of~(\ref{dfads}).
\endproof
\begin{remark} (i) Set $\eta_1=\min\{1,\eta\}$. Suppose that for every $i$ the function $\varphi_i(t)=t^{1-\theta_i}/(1-\theta_i)$, where $\theta_i\in[0,1)$, and define $\theta=\max_{i\in I(\bar{x})}\theta_i$. Then the function $\varphi(t)$ reduces to
\[\varphi(t)=\int_0^t\max_{i\in I(\bar{x})}s^{-\theta_i}ds=\int_0^t s^{-\theta}ds=\frac{t^{1-\theta}}{1-\theta},\forall t\in(0,\eta_1),\]
where the second equality holds because $\max_{i\in I(\bar{x})}s^{-\theta_i}=s^{-\theta}$ for $s\in(0,1)$. In this way, we recovered a result by Li and Pong \cite[Theorem 3.1]{li2018calculus}, where they proved that $f(x)=\min_{1\leq i\leq m} f_i(x)$ admits KL exponent $\theta$ at $\bar{x}$, if the function $f_i$ has the concave KL property at $\bar{x}$ with KL exponent $\theta_i$ for each $i\in I(\bar{x})$.

(ii) Let us provide an example where the minimum of two lsc functions is continuous. Consider $f:\R\rightarrow\R$ and $g:\R\rightarrow\R$ given by
\begin{align*}
f(x)=\begin{cases}1-x^2,&\text{ if }x\neq0;\\0,&\text{ if }x=0.\end{cases}\text{ and }g(x)=\begin{cases}x^2/2,&\text{ if }x\neq1;\\0,&\text{ if }x=1.\end{cases}
\end{align*}
Then the function $h(x)=\min\{f(x),g(x)\}$ satisfies
\begin{align*}h(x)=\begin{cases}x^2/2,&\text{ if }|x|\leq\sqrt{2/3};\\1-x^2,&\text{ if }|x|>\sqrt{2/3},\end{cases}\end{align*}
which is a continuous function. However, it is worth noting that the minimum of lsc functions is usually not continuous. Therefore a more restrictive yet easier to be verified condition for Theorem~\ref{Thm:minmum of g-KL} is that $f_i$ is continuous for every $i$, in which case the continuity of $f$ becomes automatic and our conclusion follows similarly.
\end{remark}

Next we propose a separable sum rule for the generalized concave KL property. Our proof adapted the approach in~\cite[Theorem 3.3]{li2018calculus}, but employs nondifferentiable desingularizing functions with general forms.

\begin{theorem}\label{Thm:separable sum}Let $n_i\in\N,i=1,\ldots,m,$, and let $n=\sum_{i=1}^mn_i$. For each $i$, let $f_i:\R^{n_i}\rightarrow\overline{\R}$ be a proper and lsc function that is continuous on $\dom\partial f_i$. Furthermore, suppose that each $f_i$ has the generalized concave KL property at $\bar{x}_i\in\dom\partial f_i$ with respect to $U_i=\BB{\bar{x}_i}{\varepsilon_i}$ for $\varepsilon_i>0$, $\eta_i>0$ and $\varphi_i\in\Phi_{\eta_i}$. Let $\varepsilon=\min_{1\leq i\leq m}\varepsilon_i$ and let $\eta$ be a real number with $\eta<m\cdot\min_{1\leq i\leq m}\eta_i$. Define $\varphi:[0,\eta]\rightarrow\R$ by \[\varphi(t)=\int_0^t\max_{1\leq i\leq m}\left[(\varphi_i)^\prime_-\left(\frac{s}{m}\right)\right]ds,\forall t\in(0,\eta],\]
	and $\varphi(0)=0$. If the function $\varphi_i$ is strictly concave on $(0,\eta_i)$ for each $i\in\{1,\ldots,m\}$, then the separable sum $f(x)=\sum_{i=1}^mf_i(x_i)$ has the generalized concave KL property at $\bar{x}=(\bar{x}_1,\ldots,\bar{x}_m)$ with respect to $U=\BB{\bar{x}}{\varepsilon}$, $\eta$ and $\varphi$. We take $\eta=\infty$ if $\eta_i=\infty$ for every $i$.
\end{theorem}
\proof  For each $i$, the function $\varphi_i\in\Phi_{\eta_i}$ is finite and right-continuous at $0$ according to the assumption. Invoking Lemma~\ref{Lemma: integral of maxmum of left-d is well-defined}, one concludes that $\varphi$ is finite and right-continuous at $0$. Moreover, Lemma~\ref{Lemma:integration of increasing function is convex} implies that $\varphi$ belongs to $\Phi_\eta$ with \begin{equation}\label{dfadsfas}\varphi_-^\prime(t)\geq\max_{1\leq i\leq m}(\varphi_i)_-^\prime\left(\frac{t}{m}\right)\geq(\varphi_i)_-^\prime\left(\frac{t}{m}\right),\forall 1\leq i\leq m. \end{equation}

Define $\phi_i(t)=(\varphi_i)^\prime_-(t)$ for each $i\in\{1,\ldots,m\}$. By Lemma~\ref{left derivative inequality}, the function $\phi_i$ is strictly decreasing for each $i$ hence invertible. According to the assumption that $f_i$ has the generalized concave KL property at $\bar{x}_i$, we have for $x_i\in U_i\cap[0<f_i-f_i(\bar{x}_i)<\eta_i]$,
	\begin{equation}\label{range}\phi_i(f_i(x_i)-f_i(\bar{x}_i))\geq\frac{1}{\dist(0,\partial f_i(x_i))},\end{equation}
hence $1/\dist(0,\partial f_i(x_i))\in\dom\phi_i^{-1}$. Furthermore, since $\phi_i^{-1}$ is decreasing, we have
	\[f_i(x_i)-f_i(\bar{x}_i)=\phi_i^{-1}\left(\phi_i(f_i(x_i)-f_i(\bar{x}_i) )\right)\leq\phi_i^{-1}\left(\frac{1}{\dist(0,\partial f_i(x_i))}\right). \]
	Shrinking $\varepsilon_i$ if necessary, we assume $f_i(x_i)<f_i(\bar{x}_i)+\eta_i$ whenever $x_i\in U_i$.  Therefore for every $x_i\in U_i$, we have
	\begin{align}\label{separable rule_1}
	f_i(x_i)-f_i(\bar{x}_i)\leq\phi_i^{-1}\left(\frac{1}{\dist(0,\partial f_i(x_i))}\right).
	\end{align}
In particular, the above inequality holds trivially when $f_i(x_i)\leq f_i(\bar{x}_i)$, because the right hand side is always positive.
	
Take $x=(x_1,\ldots,x_m)\in U\cap[0<f-f(\bar{x})<\eta]$ and denote by $i^*=i(x)$ the index such that \[\phi_{i^*}^{-1}\left(\frac{1}{\dist(0,\partial f_{i^*}(x_{i^*}))}\right)\geq\phi_i^{-1}\left(\frac{1}{\dist(0,\partial f_i(x_i))}\right),\forall i\in\{1,\ldots,m\}.\]
	Note that $\norm{x-\bar{x}}^2=\sum_{i=1}^m\norm{x_i-\bar{x}_i}^2$. Then $\norm{x_i-\bar{x}_i}\leq\norm{x-\bar{x}}<\varepsilon\leq\varepsilon_i$ for each $i$. For simplicity, set $r=r(x)=1/\dist(0,\partial f_{i^*}(x_{i^*}))$, where the value $r$ depends on $x$ because the index $i^*$ does. Summing (\ref{separable rule_1}) from $i=1$ to $m$ yields
	\begin{align} f(x)-f(\bar{x})&=\sum_{i=1}^m\left[f_i(x_i)-f_i(\bar{x}_i)\right]\leq\sum_{i=1}^m\phi_i^{-1}\left(\frac{1}{\dist(0,\partial f_i(x_i))}\right)\leq m\cdot\phi_{i^*}^{-1}\left(r\right). \label{separable sum_2}
	\end{align}
Let $u\in\partial f(x)$ be such that $\norm{u}=\dist(0,\partial f(x))$. Recall from Fact~\ref{Fact:subd for separable sum} that
	\[\partial f(x)=\prod_{i=1}^{m}\partial f_i(x_i).\]
Then there exists $u_i\in\partial f_i(x_i)$ for each $i$ such that $u=(u_1,\ldots,u_m)$ and consequently for every $i$, one has, $$\dist(0,\partial f(x))=\norm{u}\geq\norm{u_i}\geq\dist(0,\partial f_i(x_i)),$$
which implies that
	\begin{equation}\label{qAq}
	r\cdot\dist(0,\partial f(x))\geq1.
	\end{equation}	

Finally, we show that $\varphi$ is a desingularizing function of $f$ at $\bar{x}$. Take $x\in\BB{\bar{x}}{\varepsilon}\cap[0<f-f(\bar{x})<\eta]$. Note that the range of $\phi_{i^*}^{-1}$ satisfies $\text{ran }\phi_{i^*}^{-1}=\dom\phi_{i^*}=\dom(\varphi_{i^*})_-^\prime=(0,\eta_{i^*})$. Therefore $m\cdot\phi_{i^*}^{-1}(r)<m\cdot\eta_{i^*}$. On the other hand, recall that $\eta$ is defined to be a real number satisfying $\eta<m\cdot\min_i\eta_i\leq m\cdot \eta_{i^*}$. Therefore we need to consider two cases:
	
\textbf{Case 1}: If $m\cdot\phi_{i^*}^{-1}\left(r\right)<\eta$, then one has the following from (\ref{separable sum_2}) and the fact that $\varphi_-^\prime$ is decreasing
		\begin{align}\label{dfsfaksdals}
		\varphi_-^\prime(f(x)-f(\bar{x}))\geq\varphi_-^\prime\left(m\cdot\phi_{i^*}^{-1}\left(r\right)\right).
		\end{align}
Hence we have
	\begin{align*}
	&~~~~\varphi_-^\prime(f(x)-f(\bar{x}))\cdot\dist(0,\partial f(x))\geq\varphi_-^\prime\left(m\cdot\phi_{i^*}^{-1}\left(r \right)\right)\cdot \dist(0,\partial f(x))\\
	&\geq(\varphi_{i^*})_-^\prime\left(\phi_{i^*}^{-1}\left(r \right)\right)\cdot \dist(0,\partial f(x))=\phi_{i^*}\left(\phi_{i^*}^{-1}\left(r\right) \right)\cdot \dist(0,\partial f(x))\\
	&=r\cdot\dist(0,\partial f(x))\geq1.
	\end{align*}
	where the second inequality is implied by (\ref{dfadsfas}) and the last one holds because of (\ref{qAq}).
	
	\textbf{Case 2:} Now we consider the case where $m\cdot\phi_{i^*}^{-1}\left(r\right)\geq\eta$.  Note that $\eta<m\cdot\min\eta_i\leq m\cdot \eta_{i^*}$. Then we have $\eta/m\in\dom\phi_{i^*}$ and $\phi_{i^*}(\eta/m)<\infty$. On the other hand, the assumption
	 $m\cdot\phi_{i^*}^{-1}\left(r\right)\geq\eta$ implies that $r\leq\phi_{i^*}(\frac{\eta}{m})$. Altogether, one concludes that
\begin{align*}
&~~~~\varphi_-^\prime(f(x)-f(\bar{x}))\cdot\dist(0,\partial f(x))\geq\varphi_-^\prime(\eta)\cdot \dist(0,\partial f(x))\\
&\geq(\varphi_{i^*})_-^\prime\left(\frac{\eta}{m}\right)\cdot\dist(0,\partial f(x))\geq r\cdot\dist(0,\partial f(x))\geq1,
\end{align*}
where the second inequality is implied by (\ref{dfadsfas}) and the last one holds because of (\ref{qAq}). It is worth noting that one can take $\eta=\infty$ if $\eta_i=\infty$ for every $i$, in which case $m\cdot\phi_{i^*}^{-1}(r)<\eta$ is trivially true. Then the desired result readily follows from Case 1.\endproof
\begin{remark}\label{Rem:separable sum}  Setting $\eta<m$ and $\varphi_i(t)=t^{1-\theta_i}/(1-\theta_i)$, where $\theta_i\in(0,1)$ for every $i$, we have for $t\in(0,\eta]\subseteq(0,m)$
\begin{align*}
\varphi(t)=\int_0^t\max_{1\leq i\leq m}\left(\frac{s}{m}\right)^{-\theta_i}ds=\int_0^t\left(\frac{s}{m}\right)^{-\theta}ds=\frac{m^\theta}{1-\theta}\cdot t^{1-\theta},
\end{align*}
where $\theta=\max_{1\leq i\leq m}\theta_i$, from which a result by Li and Pong \cite[Theorem 3.3]{li2018calculus} is recovered. Note that the second equality holds because $\max_{1\leq i\leq m}t^{-\theta_i}=t^{-\theta}$ for $t\in(0,1)$ and $s/m<1$.
\end{remark}

To obtain a composition rule for the generalized concave KL property, the following technical lemma helps.
\begin{lemma}\label{Lemma: for composition rule} Let $F:\Rn\rightarrow\R^m$ be a smooth map and let $\bar{x}\in\Rn$. If the Jacobian $\nabla F(\bar{x})$ has rank $m$, then there exist $\alpha>0$ and $\varepsilon>0$ such that for $x\in\BB{\bar{x}}{\varepsilon}$,
\begin{align}\label{Formula:scalar for composition rule}
\norm{y}\leq\alpha\norm{\left(\nabla F(x) \right)^*(y)},\forall y\in\R^m\backslash\{0\}.
\end{align}
\end{lemma}
\proof By the assumption, $\nabla F(\bar{x}):\Rn\rightarrow\R^m$ is surjective continuous linear map. Then by using the open mapping theorem, there exists $\alpha>0$ such that
	\begin{equation}\label{blablabla}
	\Ball_{\R^m}\subseteq\nabla F(\bar{x})\left(\frac{\alpha}{2}\cdot\Ball_{\R^n}\right),
	\end{equation}	
where $\Ball_{\R^m}$ and $\Ball_{\Rn}$ denote the Euclidean unit balls in $\R^m$ and $\Rn$, respectively. Note that the map $F$ is assumed to be smooth. Then there exists $\varepsilon>0$ such that \[\norm{x-\bar{x}}<\varepsilon\Rightarrow\norm{\nabla F(x)-\nabla F(\bar{x})}<\frac{1}{\alpha}.\]

	We claim that for $x$ with $\norm{x-\bar{x}}<\varepsilon$, \begin{equation}\label{Formula: Ball inclusioin}
	\Ball_{\R^m}\subseteq\nabla F(x)(\alpha\Ball_{\R^n}).
	\end{equation}
	Take $u\in\Ball_{\R^m}$. We will prove the above inclusion by constructing $v\in \alpha\Ball_{\Rn}$ such that $u=\nabla F(x)(v)$. By using the inclusion (\ref{blablabla}), there exists some $v_1\in \frac{\alpha}{2}\Ball_{\R^n}$ such that $\nabla F(\bar{x})(v_1)=u$. Then we have
	\[\norm{\nabla F(x)(v_1)-u}=\norm{\nabla F(x)(v_1)-\nabla F(\bar{x})(v_1)}\leq\frac{1}{\alpha}\cdot \frac{\alpha}{2}=\frac{1}{2}.\]
	Hence $u-\nabla F(x)(v_1)\in\frac{1}{2}\Ball_{\R^m}$ and again by (\ref{blablabla}) there exists $v_2\in\frac{\alpha}{4}\Ball_{\R^n}$ such that $\nabla F(\bar{x})(v_2)=u-\nabla F(x)(v_1)$, which further implies
	\begin{align*}
	\norm{\nabla F(x)(v_1+v_2)-u}&=\norm{\nabla F(x)(v_2)-(u-\nabla F(x)(v_1))}\\
	&=\norm{\nabla F(x)(v_2)-\nabla F(\bar{x})(v_2)}\leq\frac{1}{4}.
	\end{align*}
	Repeating the same process one obtains a sequence $(v_l)_{l\in\N}$ satisfying
	\[\norm{\nabla F(x)\left(\sum_{k=1}^lv_k\right)-u}\leq\frac{1}{2^{l}}\text{ and }\norm{v_l}\leq\frac{\alpha}{2^{l}},\forall l\in\N.\]
	The latter inequality implies there exists $v\in\R^n$ such that $\lim_{l\rightarrow\infty}\sum_{k=1}^lv_k=v$ and $\norm{v}\leq\alpha$, while from the former one we have
	$\nabla F(x)(v)=u$, which proves our claim. Let $u\in\Ball_{\R^m}$ obey $\norm{u}=1$ and suppose that $u=\nabla F(x)(v)$ for some $v\in\alpha\Ball_{\R^n}$. Then we have
	
	\begin{align*}
\frac{\norm{u}}{\alpha^2}&=\frac{\ip{u}{u}}{\alpha^2}=\frac{\ip{u}{\nabla F(x) (v)}}{\alpha^2}=\frac{\ip{(\nabla F(x))^*(u)}{v}}{\alpha^2}\\
&\leq\frac{\norm{(\nabla F(x))^*(u)}\norm{v}}{\alpha^2}\leq\frac{\alpha\norm{(\nabla F(x))^*(u)}}{\alpha^2}\\
&=\frac{\norm{(\nabla F(x))^*(u)}}{\alpha},
	\end{align*}
	which implies that $\norm{u}\leq \alpha\norm{\left(\nabla F(x)\right)^*(u)}$.
Let $u=y/\norm{y}$ for nonzero $y\in\R^m$. Then one has $\norm{y}\leq\alpha\norm{\left(\nabla F(x)\right)^*(y)}$, as claimed.
\endproof

\begin{theorem}\label{Thm: Composition rule}
	Suppose that $f(x)=g(F(x))$, where $g:\R^m\rightarrow\overline{\R}$ is proper lsc and $F:\R^n\rightarrow\R^m$ is a smooth map. Let $\bar{x}\in\dom\partial f$ and let $\varepsilon_0>0$. Assume that $g$ has the generalized concave KL property at $F(\bar{x})$ with respect to $U_0=\BB{F(\bar{x})}{\varepsilon_0}$, $\eta>0$ and $\varphi\in\Phi_\eta$. If the Jacobian $\nabla F(\bar{x})$ has rank $m$, then there exist $\alpha>0$ and $\varepsilon_1>0$ such that (\ref{Formula:scalar for composition rule}) holds. Furthermore, there exits $\varepsilon\in(0,\min\{\varepsilon_0,\varepsilon_1\}]$ such that $f$ has the generalized concave KL property at $\bar{x}$ with respect to $U_1=\BB{\bar{x}}{\varepsilon}$, $\eta>0$ and $\alpha\cdot\varphi\in\Phi_\eta$.
\end{theorem}
\proof 	By assumption for $y\in\BB{F(\bar{x})}{\varepsilon_0}\cap[0<g-g(F(\bar{x}))<\eta]$ we have
	\begin{align}\label{Formula: composition rule-1}
	\varphi_-^\prime(g(y)-g(F(\bar{x})))\cdot\dist(0,\partial g(y))\geq1.
	\end{align}
	On the other hand, Lemma \ref{Lemma: for composition rule} implies that there exist $\alpha>0$ and $\varepsilon_1>0$ such that for $x\in\BB{\bar{x}}{\varepsilon_1}$,
		\begin{align}\label{dddd}
		\norm{y}\leq\alpha\norm{\left(\nabla F(x) \right)^*(y)},\forall y\in\R^m\backslash\{0\}.
		\end{align}
	Moreover, by applying Fact~\ref{Fact: chain rule}, one has $\partial f(x)=\nabla F(x)^*\partial g(F(x))=\{\nabla F(x)^*u:u\in\partial g(F(x))\}$. Let $v\in\partial f(x)$ be such that $\norm{v}=\dist(0,\partial f(x))$. Then we have for some $u\in\partial g(F(x))$
		\begin{equation}\label{Formula:sss}
		\dist(0,\partial f(x))=\norm{v}=\norm{\left(\nabla F(x)\right)^*u}\geq\frac{\norm{u}}{\alpha}\geq\frac{\dist(0,\partial g(F(x)))}{\alpha},
		\end{equation}
	where the first inequality is implied by (\ref{dddd}). Suppose that $\norm{F(x)-F(\bar{x})}<\varepsilon_0$ whenever $\norm{x-\bar{x}}<\varepsilon_2$ for some $\varepsilon_2$ and set $\varepsilon=\min\{\varepsilon_0,\varepsilon_1,\varepsilon_2\}$. Then for $x\in\BB{\bar{x}}{\varepsilon}\cap[0<f-f(\bar{x})<\eta]$
	\begin{align*}
	&~~~~(\alpha\cdot\varphi)_-^\prime(f(x)-f(\bar{x}))\dist(0,\partial f(x))\\
	&\geq\varphi^\prime_-(g(F(x))-g(F(\bar{x}))\dist(0,\partial g(F(x)))\geq1,
	\end{align*}
		where the last inequality follows from (\ref{Formula: composition rule-1}) and (\ref{Formula:sss}).
\endproof
\begin{remark} Theorem~\ref{Thm: Composition rule} generalizes~\cite[Theorem 3.2]{li2018calculus}, where Li and Pong obtained a similar composition rule with an additional assumption that desingularizing functions have the form $\varphi(t)=c\cdot t^{1-\theta}$ for $c>0$ and $\theta\in[0,1)$. Moreover, they proved~(\ref{Formula: Ball inclusioin}) by applying the Lyusternik-Graves theorem \cite[Theorem 1.57]{mor2018variational} to the continuous map $x\mapsto\nabla f(x)$. In contrast, we provided a different proof by using Lemma~\ref{Lemma: for composition rule}. 
\end{remark}

\begin{corollary}\label{Cor: speical composition} Suppose that $f(x)=g(Ax-b)$, where $A\in\R^{m\times n}$ has rank $m$ and $b\in\R^m$. Let $\bar{x}\in\dom\partial f$ and $\varepsilon>0$. Set $r=\sqrt{\lambda_{\min}(AA^*)}>0$. If $g$ has the generalized concave KL property at $A\bar{x}-b$ with respect to $U_1=\BB{A\bar{x}-b}{\varepsilon}$, $\eta>0$ and $\varphi(t)\in\Phi_\eta$, then $f$ has the generalized concave KL property at $\bar{x}$ with respect to $U_2=\BB{\bar{x}}{\varepsilon/\norm{A}}$, $\eta>0$ and $\varphi(t)/r$. Note that $U_2=\Rn$ if $U_1=\R^m$.
\end{corollary}
\proof Notice that $A$ is surjective. Hence~\cite[Exercise 1.53]{mor2018variational} implies that
\begin{align*}
\norm{A^*y}\geq r\norm{y},
\end{align*}
which means that $F(x)=Ax-b$ satisfies (\ref{Formula:scalar for composition rule}) for every $x\in\Rn$ with $\alpha=1/r$ and $\varepsilon_1=\infty$. Then applying a similar argument in Theorem \ref{Thm: Composition rule} completes the proof. \endproof

Finally, we highlights our results with two examples.
The first example shows that the Li-Pong calculus rule for minimum of functions~\cite[Theorem 3.1]{li2018calculus} fails, while Theorem~\ref{Thm:minmum of g-KL} is still applicable.
\begin{example}\label{ex:min rule fails} Set $\tilde{\varphi}(t)=\sqrt{-1/\ln(t)}$ for $t>0$ and $\tilde{\varphi}(0)=0$. Define $f(x)=\exp(-1/x^2)$ for $x\neq0$ and $f(0)=0$. Let $g(x)=|x|$ and $h(x)=\min\{f(x),g(x)\}$. Then there exists $\varepsilon>0$ such that $h$ has the generalized concave KL property at $\bar{x}=0$ with respect to $U=(-\varepsilon,\varepsilon)$, $\eta=\exp(-3/2)$ and $\varphi(t)=\tilde{\varphi}(t)$.
\end{example}
\proof Clearly $g$ has the generalized concave KL property at $\bar{x}$ with respect to $U=\R$, $\eta=\infty$ and $\varphi(t)=t$. Proposition~\ref{Prop: the usual form doesn't always work}(i) states that $f$ has the generalized concave KL property with respect to $U=(-\sqrt{2/3},\sqrt{2/3})$, $\eta=\exp(-3/2)$ and $\tilde{\varphi}(t)$. The desired result follows immediately from Theorem~\ref{Thm:minmum of g-KL}. In particular, we have
\[\varphi(t)=\int_0^t\max\{(\tilde{\varphi})^\prime(s),1\}ds=\tilde{\varphi}(t),\]
where the last equality holds because $\tilde{\varphi}^\prime(t)=1/(2t\sqrt{(-\ln(t))^3})>1$ for $t\in(0,\eta)$. \endproof
\begin{remark} Recall that calculus rules in~\cite{li2018calculus} assume desingularizing functions to have the special form $\varphi(t)=c\cdot t^{1-\theta}$ for $c>0$ and $\theta\in[0,1)$. By using Proposition~\ref{Prop: the usual form doesn't always work}(ii), it is easy to see that~\cite[Theorem 3.1]{li2018calculus} fails.
\end{remark}

The next example shows that applying generalized calculus rules to the exact modulus
(recall Definition~\ref{Def: g-KL}) could lead to a smaller desingularizing function.

\begin{example}\label{ex:toy example} Define $f(x_1,x_2)=-\ln(1-x_1^2)-\ln(1-x_2^2)$ for $(x_1,x_2)\in(-1,1)\times(-1,1)$. Let $\bar{x}=0$. Then $f$ has the generalized concave KL property at $\bar{x}$ with respect to $U=\BB{\bar{x}}{1}$, $\eta=\infty$ and $\varphi(t)=2\sqrt{1-\exp(-t/2)}$.\end{example}
\proof Set $f_i(x_i)=-\ln(1-x_i^2)$ for $i=1,2$. Then simple calculation yields that $f_i$ has the generalized concave KL property at $0$ with respect to $U_i=(-1,1)$, $\eta_i=\infty$ and the exact modulus $\varphi_i(t)=\sqrt{1-\exp(-t)}$. Note that $\varphi_i$ is strictly concave for every $i$. Then, invoking Theorem~\ref{Thm:separable sum} with $\varepsilon_i=1$ and $\eta_i=\infty$, one concludes that $f$ has the generalized concave KL property at $0$ with respect to $U=\BB{\bar{x}}{1}$, $\eta=\infty$ and $\varphi(t)=2\sqrt{1-\exp(-t/2)}$.\endproof
\begin{remark}\label{rem: smaller desingularizing function}Since $f_i^{\prime\prime}(x_i)=(2+2x_i^2)/(1-x_i^2)^2\geq2$, where $f_i$ is given in the proof of Example~\ref{ex:toy example}, the separable sum $f$ is $2$-strongly convex, i.e., $f-\norm{\cdot}^2$ is convex. By applying~\cite[Example 6]{Bolte2014}, $f$ has the generalized concave KL property $\bar{x}=0$ with respect to $U=\dom f$, $\eta=\infty$ and $\varphi_1(t)=2\sqrt{t}$. It is easy to see that $\varphi(t)\leq2\sqrt{t/2}\leq\varphi_1(t)$.
\end{remark}

\section{Conclusion}\label{Section:Conclusion}

We established several calculus rules for the generalized concave KL property, which generalize Li and Pong's calculus rules for the KL exponent~\cite{li2018calculus}; see also its early preprint~\cite{li2016calculus}. Compared to their results, ours do not assume desingularizing functions to have any specific forms nor differentiable. Such generalization is motivated by the recent discovery that the optimal concave desingularizing function has various forms and may not be differentiable. Let us end this paper by discussing some directions for future work:
\begin{itemize}
	\item It is tempting to estimate the exact modulus for concrete optimization models using tools developed in this paper.
	\item The classic approach to determine convergence rates of algorithms under the KL assumption is to find an appropriate KL exponent $\theta$, given that desingularizing functions have the form $\varphi(t)=c\cdot t^{1-\theta}$ for $c>0$ and $\theta\in[0,1)$. Now that we have found a way to calculate smaller desingularizing functions with possibly different forms, it is interesting to explore whether this framework leads to a new and possibly sharper analysis of convergence rate.
\end{itemize}

\section*{Acknowledgments}
XW and ZW were partially supported by NSERC Discovery Grants. The authors thank Dr. Heinz H. Bauschke for many help discussions.
\appendix
\addcontentsline{toc}{section}{Appendices}
\renewcommand{\thesubsection}{\Alph{subsection}}
\section*{Appendix}
\counterwithin{theorem}{subsection}

\subsection{Supplementary lemmas}

\begin{lemma}\label{Lemma: integral of maxmum of left-d is well-defined} Let $m\in\N$ obey $m\geq2$. For each $i\in\{1,\ldots,m\}$, let $h_i:(0,\infty)\rightarrow\R_+$ be such that $\lim_{s\rightarrow0^+}h(s)=\infty$. Suppose that for each $i$ the function $\varphi_i:[0,\infty)\rightarrow\R_+$ given by $\varphi_i(t)=\int_0^th_i(s)ds$ for $t\in(0,\infty)$ and $\varphi_i(0)=0$, is finite and right-continuous at $0$. Then the function $\varphi:[0,\infty)\rightarrow\R_+$, \[(\forall 0<t<\eta)~t\mapsto\int_0^t\max_{1\leq i\leq m}h_i(s)ds,\] and $\varphi(0)=0$, is finite and right-continuous at $0$.
\end{lemma}
\proof Note that $\varphi$ is an improper integral. Hence we have for $t>0$ \[\varphi(t)=\lim_{u\rightarrow0^+} \int_u^t\max_{1\leq i\leq m}h_i(s)ds.\]
Let $m=2$. Then by using the inequality $\max\{\alpha,\beta\}\leq\alpha+\beta$ for $\alpha,\beta\geq0$, one has
\begin{align*}
	\max\{h_1(s),h_2(s)\}\leq h_1(s)+h_2(s).
\end{align*}
Hence we have for $t>0$
\begin{align*}
	\varphi(t)&\leq\lim_{u\rightarrow0^+}\int_u^th_1(s)+h_2(s)ds\\
	&=\varphi_1(t)+\varphi_2(t)-\lim_{u\rightarrow0^+}\left[\varphi_1(u)+\varphi_2(u)\right]\\
	&=\varphi_1(t)+\varphi_2(t)<\infty,
\end{align*}
where the last equality is implied by the right-continuity of $\varphi_1$ and $\varphi_2$ at $0$. Taking $t\rightarrow0^+$, one gets $\lim_{t\rightarrow0^+}\varphi(t)=0$. The desired result then follows from a simple induction.\endproof

\begin{lemma}\label{left derivative inequality} Let $\eta\in(0,\infty]$ and let $\varphi\in\Phi_\eta$. Then the following
hold:
	\begin{itemize}
		\item[(i)] Let $t>0$. Then $\varphi(t)=\lim_{u\rightarrow0^+}\int_u^t\varphi_-^\prime(s)ds=\int_0^t\varphi_-^\prime(s)ds$.
		\item [(ii)] The function $t\mapsto\varphi_-^\prime(t)$ is decreasing and $\varphi_-^\prime(t)>0$ for $t\in(0,\eta)$. If in addition $\varphi(t)$ is strictly concave, then $t\mapsto\varphi_-^\prime(t)$ is strictly decreasing.
		\item[(iii)] For $0\leq s<t<\eta$, $\varphi^\prime_-(t)\leq\frac{\varphi(t)-\varphi(s)}{t-s}$.
	\end{itemize}
\end{lemma}
\proof (i) Invoking Lemma~\ref{Lemma: properties of convex functions on the line}(ii) yields
\[\varphi(t)=\lim_{u\rightarrow0^+}\big(\varphi(t)-\varphi(u)\big)=\lim_{u\rightarrow0^+}\int_{u}^t\varphi_-^\prime(s)ds<\infty,\]
where the first equality holds because $\varphi$ is right-continuous at $0$ with $\varphi(0)=0$. Let $(u_n)_{n\in\N}$ be a decreasing sequence with $u_1<t$ such that $u_n\rightarrow0^+$ as $n\rightarrow\infty$. For each $n$, define $h_n:(0,t]\rightarrow\R_+$ by $h_n(s)=\varphi_-^\prime(s)$ if $s\in(u_n,t]$ and $h_n(s)=0$ otherwise. Then the sequence $(h_n)_{n\in\N}$ satisfies: (a) $h_n\leq h_{n+1}$ for every $n\in\N$; (b) $h_n(s)\rightarrow\varphi_-^\prime(s)$ pointwise on $(0,t)$; (c) The integral $\int_0^th_n(s)ds=\int_{u_n}^t\varphi_-^\prime(s)ds=\varphi(t)-\varphi(u_n)\leq \varphi(t)-\varphi(0)<\infty$ for every $n\in\N$. Hence the monotone convergence theorem implies that
\[\lim_{u\rightarrow0^+}\int_{u}^t\varphi_-^\prime(s)ds=\lim_{n\rightarrow\infty}\int_{u_n}^t\varphi_-^\prime(s)ds=\lim_{n\rightarrow\infty}\int_0^th_n(s)ds=\int_0^t\varphi_-^\prime(s)ds.\]

(ii) According to Lemma~\ref{Lemma: properties of convex functions on the line}(i), the function $t\mapsto\varphi_-^\prime(t)$ is decreasing. Suppose that $\varphi_-^\prime(t_0)=0$ for some $t_0\in(0,\eta)$. Then by the monotonicity of $\varphi_-^\prime$ and (i), we would have $\varphi(t)-\varphi(t_0)=\int_{t_0}^t\varphi_-^\prime(s)ds\leq(t-t_0)\varphi_-^\prime(t_0)=0$ for $t>t_0$, which contradicts to the assumption that $\varphi$ is strictly increasing.

(iii) For $0<s<t<\eta$, applying Lemma~\ref{Lemma: properties of convex functions on the line}(iii) to the convex function $-\varphi$ yields that $-\varphi(s)+\varphi(t)\geq-\varphi_-^\prime(t)(s-t)\Leftrightarrow\varphi_-^\prime(t)\leq\big(\varphi(t)-\varphi(s)\big)/(t-s)$. The desired inequality then follows from the right-continuity of $\varphi$ at $0$.\endproof

\begin{lemma}\label{Lemma:integration of increasing function is convex}
	Let $\eta\in(0,\infty]$ and let $h:(0,\eta)\rightarrow\R_+$ be a positive-valued decreasing function. Define $\varphi(t)=\int_0^th(s)ds$ for $t\in(0,\eta)$ and set $\varphi(0)=0$. Suppose that $\varphi(t)<\infty$ for $t\in(0,\eta)$. Then $\varphi$ is a strictly increasing concave function on $[0,\eta)$ with $$\varphi_-^\prime(t)\geq h(t)$$ for $t\in(0,\eta)$, and right-continuous at $0$. If in addition $h$ is a continuous function, then $\varphi$ is $C^1$ on $(0,\eta)$.
\end{lemma}
\proof Let $0<t_0<t_1<\eta$. Then $\varphi(t_1)-\varphi(t_0)=\int_{t_0}^{t_1}h(s)ds\geq(t_1-t_0)\cdot h(t_1)>0$, which means $\varphi$ is strictly increasing. Applying~\cite[Theorem 6.79]{Stromberg}, one concludes that $\varphi(t)\rightarrow\varphi(0)=0$ as $t\rightarrow0^+$.\endproof

\subsection{Properties of the generalized concave KL property and its associated exact modulus}\label{AppendixA}

In this subsection, we recall some pleasant properties of the generalized concave KL property and its associated exact modulus. Results in this subsection can be found in~\cite{wang2020}. Here we provide detailed proofs for the sake of self-containedness.

\begin{proposition}\label{Prop: g-KL modulus is optimal} Let $f:\Rn\rightarrow\overline{\R}$ be proper lsc and let $\bar{x}\in\dom\partial f$. Let $U$ be a nonempty neighborhood of $\bar{x}$ and $\eta\in(0,\infty]$. Let $\varphi\in\Phi_\eta$ and suppose that $f$ has the generalized concave KL property at $\bar{x}$ with respect to $U$, $\eta$ and $\varphi$. Then the exact modulus of the generalized concave KL property of $f$ at $\bar{x}$ with respect to $U$ and $\eta$, denoted by $\tilde{\varphi}$, is well-defined and satisfies \[\tilde{\varphi}(t)\leq\varphi(t),~\forall t\in[0,\eta).\]
	Moreover, the function $f$ has the generalized concave KL property at $\bar{x}$ with respect to $U$, $\eta$ and $\tilde{\varphi}$. Furthermore, the exact modulus $\tilde{\varphi}$ satisfies
	\begin{equation*}
		\tilde{\varphi}=\inf\big\{\varphi\in\Phi_\eta:\text{$\varphi$ is a concave desingularizing function of $f$ at $\bar{x}$ with respect to $U$ and $\eta$}\big\}.
	\end{equation*}
\end{proposition}
\proof  Let us show first that $\tilde{\varphi}(t)\leq\varphi(t)$ on $[0,\eta)$, from which the well-definedness of $\tilde{\varphi}$ readily follows. If $U\cap[0<f-f(\bar{x})<\eta]=\emptyset$, then by our convention $\tilde{\varphi}(t)=0\leq\varphi(t)$ for every $t\in[0,\eta)$. Therefore we proceed with assuming $U\cap[0<f-f(\bar{x})<\eta]\neq\emptyset$. By assumption, one has for $x\in U\cap[0<f-f(\bar{x})<\eta]$,
\[\varphi_-^\prime\big(f(x)-f(\bar{x})\big)\cdot\dist\big(0,\partial f(x)\big)\geq1.\]
which guarantees that $\dist\big(0,\partial f(x)\big)>0$. Fix $s\in(0,\eta)$ and recall from Lemma~\ref{left derivative inequality}(ii) that $\varphi_-^\prime(t)$ is decreasing. Then for $x\in  U\cap[0<f-f(\bar{x})<\eta]$ with $s\leq f(x)-f(\bar{x})$ we have \[\dist^{-1}\big(0,\partial f(x)\big)\leq\varphi_-^\prime\big(f(x)-f(\bar{x})\big)\leq\varphi_-^\prime(s).\]
Taking the supremum over all $x\in U\cap[0<f-f(\bar{x})<\eta]$ satisfying $s\leq f(x)-f(\bar{x})$ yields
\[h(s)\leq\varphi_-^\prime(s),\]
where $h(s)=\sup\big\{\dist^{-1}\big(0,\partial f(x)\big): x\in U\cap[0<f-f(\bar{x})<\eta],s\leq f(x)-f(\bar{x})\big\}$. If $\lim_{s\rightarrow0^+}h(s)=\infty$, then one needs to treat $\tilde{\varphi}(t)$ as an improper integral. For $t\in(0,\eta)$, one has
\[\tilde{\varphi}(t)=\lim_{u\rightarrow0^+}\int_u^th(s)ds\leq\lim_{u\rightarrow0^+}\int_u^t\varphi_-^\prime(s)ds=\varphi(t)<\infty,\]
where the last equality follows from Lemma~\ref{left derivative inequality}. If $\lim_{s\rightarrow0^+}h(s)<\infty$, then the above argument still applies.

Recall that $\dist\big(0,\partial f(x)\big)>0$ for every $x\in U\cap[0<f-f(\bar{x})<\eta]$. Hence $h(s)$ is positive-valued. Take $s_1,s_2\in(0,\eta)$ with $s_1\leq s_2$. Then for $x\in U\cap[0<f-f(\bar{x})<\eta]$, one has
\[s_2\leq f(x)-f(\bar{x})\Rightarrow s_1\leq f(x)-f(\bar{x}),\]
implying that $h(s_2)\leq h(s_1)$. Therefore $h(s)$ is decreasing. Invoking Lemma~\ref{Lemma:integration of increasing function is convex}, one concludes that $\tilde{\varphi}\in\Phi_\eta$ and $\varphi_-^\prime(t)\geq h(t)$ for every $t\in(0,\eta)$.

Let $t\in(0,\eta)$. Then for $x\in  U\cap[0<f-f(\bar{x})<\eta]$ with $t=f(x)-f(\bar{x})$ we have
\[\tilde{\varphi}_-^\prime\big(f(x)-f(\bar{x})\big)\geq h(t)\geq\dist^{-1}\big(0,\partial f(x)\big),\]
where the last inequality is implied by the definition of $h(s)$, from which the generalized concave KL property readily follows because $t$ is arbitrary.

Recall that $\varphi$ is an arbitrary concave desingularizing function of $f$ at $\bar{x}$ with respect to $U$ and $\eta$, and $\tilde{\varphi}(t)\leq\varphi(t)$ for all $t\in[0,\eta)$. Hence one has
\[\tilde{\varphi}\leq\inf\big\{\varphi\in\Phi_\eta:\text{$\varphi$ is a concave desingularizing function of $f$ at $\bar{x}$ with respect to $U$ and $\eta$}\big\}.\]
On the other hand, the converse inequality holds as $\tilde{\varphi}$ is a concave
desingularizing function of $f$ at $\bar{x}$ with respect to $U$ and $\eta$.\endproof

The exact modulus is not necessarily differentiable and may have different forms depending on the function of interest.
\begin{example}\label{ex: nondifferentiable modulus} 
Let $\rho>0$. Consider the function given by
		\begin{align*}
			f(x)=
			\begin{cases}
				2\rho|x|-3\rho^2/2,&\text{if }|x|>\rho;\\
				|x|^2/2, &\text{if }|x|\leq\rho.
			\end{cases}
		\end{align*}
		Then the function
		\begin{align*}
			\tilde{\varphi}(t)=\begin{cases}\sqrt{2t},&\text{if }0\leq t\leq \rho^2/2;\\ t/(2\rho)+3\rho/4, &\text{if }t>\rho^2/2,\end{cases}
		\end{align*}
		is the exact modulus of the generalized concave KL property of $f$ at $\bar{x}=0$ with respect to $U=\R$ and $\eta=\infty$.
\end{example}
\proof For $x\neq0$, one has
\begin{align*}
	\dist^{-1}\big(0,\partial f(x)\big)=\begin{cases}1/|x|,&\text{if }0<|x|\leq\rho;\\1/2\rho,&\text{if }|x|>\rho. \end{cases}
\end{align*}
It follows that for $s\in(0,\rho^2/2]$,
\begin{align*}
	h(s)&=\sup\big\{\dist^{-1}\big(0,\partial f(x)\big):x\in\R\cap[0<f<\infty],s\leq f(x)\big\}\\
	&=\sup\big\{\dist^{-1}\big(0,\partial f(x)\big):|x|\geq\sqrt{2s}\big\}=1/\sqrt{2s},
\end{align*}
and for $s>\rho^2/2$
\begin{align*}
	h(s)&=\sup\big\{\dist^{-1}\big(0,\partial f(x)\big):x\in\R\cap[0<f<\infty],s\leq f(x)\big\}\\
	&=\sup\big\{\dist^{-1}\big(0,\partial f(x)\big):x\neq0,|x|\geq s/(2\rho)+3\rho/4\big\}=1/(2\rho),
\end{align*}
from which the desired result readily follows.\endproof

The next proposition restates an example from~\cite[Section 1]{bolte2007lojasiewicz} in our extended framework. No proof was given in~\cite[Section 1]{bolte2007lojasiewicz}, thus we prove it here for the sake of completeness.
\begin{proposition}\label{Prop: the usual form doesn't always work} Define  $f:\R\rightarrow\R$ by
	$f(x)=e^{-1/x^2}$ for $x\neq0$ and $f(0)=0$. Then the following hold:

(i) The function $\tilde{\varphi}_2(t)=\sqrt{-1/\ln(t)}$ for $t>0$ and $\tilde{\varphi}_2(0)=0$
is the exact modulus of generalized KL property of $f$ at $\bar{x}=0$ with respect to $U=(-\frac{2}{3},\frac{2}{3})$ and $\eta=\exp(-\frac{3}{2})$.

(ii) For every $c>0$ and $\theta\in[0,1)$, the function $\varphi(t)=c\cdot t^{1-\theta}$ cannot be a desingularizing function of the generalized concave KL property of $f$ at $0$ with respect to any neighborhood $U\ni0$ and $\eta\in(0,\infty]$.
\end{proposition}
\proof (i) For $0<s\leq\exp(-3/2)$, $s\leq\exp(-1/x^2)\Leftrightarrow|x|\geq\sqrt{-1/\ln{s}}$. Thus we have
\begin{align*}
	h_2(s)&=\sup\{|f^\prime(x)|^{-1}: x\in U_2\cap[0<g-g(\bar{x})<\eta_2], s\leq g(x)-g(\bar{x})\}\\
	&=\sup\{|2x^{-3}\exp(-1/x^2)|^{-1}: \sqrt{-1/\ln{s}}\leq|x|<2/3\}=(-\ln(s))^{-3/2}/(2s).
\end{align*}
Hence $\tilde{\varphi}_2(t)=\sqrt{-1/\ln(t)}$ for $t>0$.

(ii) Suppose to the contrary that there were $c>0$ and $\theta\in[0,1)$ such that $f$ has the generalized concave KL property at $0$ with respect to some $U\ni0$ and $\eta>0$ and $\varphi(t)=c\cdot t^{1-\theta}$. Taking the intersection if necessary, assume without loss of generality that $U\cap[0<f<\eta]\subseteq(-\sqrt{2/3},\sqrt{2/3})\cap[0<f<e^{-3/2}]$. Then $f$ is convex and $C^1$ on $U\cap[0<f<\eta]$. Applying a similar argument as in~Example~\ref{ex: nondifferentiable modulus}, one concludes that the exact modulus of the generalized concave KL property of $f$ at $\bar{x}$ with respect to $U$ and $\min\{\eta,e^{-3/2}\}$ is also $\tilde{\varphi}(t)$. Hence Proposition~\ref{Prop: g-KL modulus is optimal} implies that
\begin{equation}\label{dddddd}
	\tilde{\varphi}(t)\leq\varphi(t)=c\cdot t^{1-\theta},~\forall t\in(0,\min\{\eta,e^{-3/2}\}).
\end{equation}
Let $s>0$. Then one has $s=\tilde{\varphi}(t)\Leftrightarrow t=e^{-1/s^2}$, which further implies that
\[\limsup_{t\rightarrow0^+}\frac{\tilde{\varphi}(t)}{t^{1-\theta}}=\limsup_{s\rightarrow0^+}\frac{s}{e^{-(1-\theta)/s^2}}=\limsup_{s\rightarrow0^+}\frac{e^{(1-\theta)/s^2}}{s^{-1}}=\infty,\]
which contradicts to (\ref{dddddd}).\endproof

\bibliographystyle{siam}
\bibliography{KL_modulus_reference}
\end{document}